\documentclass{birkmult}
 \newtheorem{thm}{Theorem}[section]

 \newtheorem{conj}[thm]{Conjecture}
 \newtheorem{prop}[thm]{Proposition}
 
 \theoremstyle{definition}
 \newtheorem{defn}[thm]{Definition}
 \theoremstyle{remark}

 \numberwithin{equation}{section}
\newcommand{\R}{\mathbb{R}}
\def\calC{{\mathcal C}}
\def\calT{{\mathcal T}}
\def\calL{{\mathcal L}}
\def\STAB{{\rm STAB}}
\def\SENS{{\rm SENS}}
\def\Var{{\rm Var}}
\def\E{{\rm E}}
\def\P{{\rm P}}
\newcommand{\Z}{\mathbb{Z}}
\def\QED{\qed\medskip}
\def\ew#1{w_{#1}}
\begin{document}
%
%
%
%
%
%
%
%
%
\title[A survey on dynamical percolation]
 {A survey of dynamical percolation}
\author[Jeffrey E. Steif]{Jeffrey E. Steif}
\address{%
Mathematical Sciences \\
Chalmers University of Technology \\
and \\
Mathematical Sciences \\
G\"{o}teborg University \\
SE-41296 Gothenburg, Sweden}

\email{steif@chalmers.se}

\thanks{Research partially supported by the Swedish Natural Science 
Research Council and the G\"{o}ran Gustafsson 
Foundation for Research in Natural Sciences and Medicine}


\date{October 18, 2008}

\dedicatory{Dedicated to the memory of Oded Schramm, who has been a 
great inspiration to me and with whom it has been a great honor and 
privilege to work.}

\begin{abstract}
Percolation is one of the simplest and nicest models in probability 
theory/statistical mechanics which exhibits critical phenomena.
Dynamical percolation is a model where a simple
time dynamics is added to the (ordinary) percolation model. 
This dynamical model exhibits very interesting behavior.
Our goal in this 
survey is to give an overview of the work in dynamical 
percolation that has been done (and some of 
which is in the process of being written up). 
\end{abstract}

\maketitle
\tableofcontents

\section {Introduction} 


Within the very large and beautiful world of probability theory, 
percolation theory has been an especially attractive subject, an area
in which the major problems are easily stated but whose solutions,
when they exist, often require ingenious methods.
Dynamical percolation is a model where we add a simple time
dynamics to the percolation model.  It happened to
have turned out that this new model is very rich
and much can be said about it. It is hoped that this survey
will provide a good guide to the literature from which point the reader
can study the various papers on the subject of dynamical percolation.

\subsection {The ordinary percolation model}
In the standard percolation model,
we take an infinite connected locally finite graph
$G$, fix $p\in[0,1]$ and let each edge (bond) of $G$ be, independently of
all others, open with probability $p$ and closed with probability $1-p$.
Write $\pi_p$ for this product measure. In percolation theory, one studies 
the structure of the connected components (clusters) of the random subgraph 
of $G$ consisting of all sites and all open edges.
The first question that can be asked is the existence of an 
infinite connected component (in which case we say percolation occurs). 
Writing $\calC$ for this latter event, Kolmogorov's 0-1 law tells us that
the probability of $\calC$ is, for fixed $G$ and $p$, either 0 or 1.
Since $\pi_p(\calC)$ is nondecreasing in $p$, there
exists a critical probability $p_c=p_c(G)\in[0,1]$ such that
\[
\pi_p(\calC)=\left\{
\begin{array}{ll}
0 & \mbox{for } p<p_c \\
1 & \mbox{for } p>p_c.
\end{array} \right.
\]
At $p=p_c$, we can have either $\pi_p(\calC)=0$ or $\pi_p(\calC)=1$, 
depending on $G$. Site percolation is defined analogously where
the vertices are retained independently with probability $p$ and all the
edges between retained vertices are retained.

For the $d$-dimensional cubic lattice, which is the graph classically studied,
the standard reference is
\cite{Grimmett}. For the study of percolation on general graphs,
see \cite{LyonsPeres}. For a study of critical percolation on the hexagonal
lattice, see \cite{W}.

While a tremendous amount is known about this model, we highlight some 
key facts. In the following, $\Z^d$ denotes the $d$-dimensional integers
with the usual graph structure where points at distance 1 are neighbors.

$\bullet$ 
For $\Z^2$, Harris \cite{Ha} established in 1960
that there is no percolation
at $1/2$ and Kesten \cite{Ke} established in 1980
that there is percolation for 
$p >1/2$. In particular, we see there is no percolation at the critical 
value for $\Z^2$.

$\bullet$ For $\Z^d$ with $d\ge 19$, Hara and Slade \cite{HS}
established that there is no percolation at the critical value. 
(Percolation connoisseurs
will object to this description and say that
they, in addition, essentially established this result for 
$d\ge 7$.)

$\bullet$ Aizenman, Kesten and Newman \cite{AKN} proved that on 
$\Z^d$,
when there is percolation, there is then a unique infinite cluster. 
Later Burton and Keane  \cite{BK}
established this in much greater generality but
more importantly found a much simpler proof of this theorem.

\medskip

We end this subsection with an important open question.

\medskip

$\bullet$ 
For $\Z^d$, for intermediate dimensions, such as
$d= 3$, is there percolation at the critical value?

\subsection {The dynamical percolation model}

H\"{a}ggstr\"{o}m, Peres and Steif~\cite{HPS} initiated the study of
dynamical percolation. In fact, Olle H\"{a}ggstr\"{o}m and I came up
with the model inspired by a question that Paul Malliavin
asked at a lecture I gave at the Mittag Leffler Institute in 1995.
This model was invented independently by Itai Benjamini.
In this model, with $p$ fixed, the edges of $G$ switch back
and forth according to independent 2-state
continuous time Markov chains where closed switches to open at rate $p$
and open switches to closed at rate $1-p$. Clearly, $\pi_p$ is the 
unique stationary
distribution for this Markov process. The general question studied in
dynamical percolation is whether, when we start with distribution $\pi_p$,
there exist atypical times at which the percolation
structure looks markedly different than that at a fixed time. In almost all
cases, the term  ``markedly different'' refers to the existence or
nonexistence of an infinite connected component. Dynamical percolation
on site percolation models is defined analogously.

It turns out that dynamical percolation is a very interesting model 
which exhibits quite new phenomena. It also connects up very much with
the recent and important developments in 2-d critical percolation.

In the paper \cite{Haggsurvey}, H\"{a}ggstr\"{o}m gives a description of
the early results in dynamical percolation and some open questions; most of
these earlier results will be repeated here. We also mention the
following upcoming paper \cite{Garban} by Garban which will discuss 
Oded Schramm's contributions to the area of noise sensitivity; this
paper will have some overlap with the present paper.

\subsection {Outline of paper}

In Section \ref{sect:firstresults}, 
we will start with a discussion of some of the first 
results obtained for dynamical percolation including 
(1) ``triviality away from criticality'', (2) the existence of 
graphs with ``exceptional times'' and (3) nonexistence of 
``exceptional times'' for high dimensional lattices.
In Section \ref{sect:tree}, we present the fairly refined results 
that have been obtained for trees that do not percolate at criticality 
in ordinary percolation.

We present in Section \ref{sect:NS} the basic elements of noise 
sensitivity and influence together with the Fourier spectrum which 
is the key tool for the analysis of these concepts.
These concepts are both interesting in themselves and are key in the
``second moment method'' approach to proving the existence of
exceptional times for dynamical percolation in 2 dimensions as well
as to computing the Hausdorff dimension of the set of exceptional times.
Section \ref{sect:crit} briefly reviews results concerning
critical exponents for ordinary critical percolation since these
are crucial to the study of dynamical percolation in 2 dimensions.

In Section \ref{sect:SS}, we outline the proofs from \cite{SS}
for the hexagonal lattice 
of the existence of exceptional times for dynamical percolation (which also
yields a nonsharp lower bound on the Hausdorff dimension)
and obtain a lower bound for the noise sensitivity
exponent for crossing probabilities. 
(For the square lattice, positivity of the noise sensitivity
exponent was established but not the existence of exceptional
times.) This method managed to estimate
the Fourier spectrum well enough to obtain the above two results. However,
the estimates of the spectrum that this method led to
were not sharp enough to yield the exact noise sensitivity exponent 
or the exact dimension of the set of exceptional times.
We present the results from \cite{GPS} in Section \ref{sect:GPS}. These
yield for the hexagonal lattice 
the exact dimension of the set of exceptional times for dynamical percolation
and the exact noise sensitivity exponent for crossing probabilities.
(For the square lattice, existence of exceptional times for 
dynamical percolation as well as sharp noise results in terms of
the number of pivotals is also established.)
These arguments are based on a completely different method to analyze
the Fourier spectrum which turns out to be sharp. 

Essentially all the results above dealt with graphs which do not 
percolate at criticality. 
When dealing with graphs which do percolate at criticality
and asking if there are exceptional times of {\em non}-percolation,
the structure of the problem is very different. A number of results
in this case are presented in Section \ref{sect:PSS}. 

In Section \ref{sect:IIC}, a relationship between
the incipient infinite cluster (in ordinary percolation) and the structure
of the infinite cluster at exceptional times in dynamical percolation
is presented. 
Recent work on the scaling limit of dynamical percolation on the 
hexagonal lattice is discussed in Section \ref{sect:SL}, . 
Finally, in Section \ref{sect:BS}, dynamical percolation for certain 
interacting particle systems are discussed.

\medskip\noindent
We end the introduction by mentioning that
there are a number of other papers, such as
\cite{BHPS}, \cite{JS} and \cite{J}, dealing with the existence 
of certain ``exceptional times'' for various models.  
On an abstract level, the questions we are asking
concerning exceptional behavior are equivalent
to questions concerning the polarity or nonpolarity of certain sets
for our (Markovian) dynamical percolation model; see \cite{Evans}.
Sometimes, the exceptional objects one looks for are of a more general
nature such as in \cite{BeSc} and \cite{ABP}. 

\section {Dynamical percolation: first results} 
\label{sect:firstresults}

This section will discuss a selected subset of the results that were obtained
in \cite{HPS} and \cite{PS}, which were the first two papers
on dynamical percolation. 

\subsection {Away from criticality.}

The first very elementary result is the following
and contained in \cite{HPS}. Let $\calC_t$ be the event that there is
an infinite cluster at time $t$ and $\P_p$ denote the probability measure
governing the whole process when parameter $p$ is used and the graph is
understood.

\begin{prop} \label{pr:noncrit} \cite{HPS}
 For any graph $G$ we have 
\begin{equation}  \begin{array}{ccl}
\P_p\bigl((\neg\, \calC_t) \mbox{ occurs for every } t\bigr)=1 & \mbox{ if } & p<p_c(G) 
   \\[1ex]
\P_p(\, \calC_t \, \mbox{ occurs for every  } \, t \, )=1  & \mbox{ if } & p>p_c(G) \, .
\end{array} 
\nonumber
\end{equation} 
\end{prop}

\medskip\noindent
Outline of proof: For the
first part, we choose $\delta$ sufficiently small so that the set of edges
which are open at {\em some time} during a time interval of length $\delta$
is still subcritical. Then on each ``$\delta$-interval'' there are no
exceptional times of percolation and we use countable additivity. The second
part is proved in the same way. $\QED$

This result suggests that the only interesting parameter is the critical
one and this is essentially true. However, before restricting ourselves in the 
continuation to the critical case, we do mention two results for the 
supercritical case, both from \cite{PS}.

The first result is a version for dynamical percolation
of the uniqueness of the infinite cluster. 

\begin{thm} \label{th:lattice}\cite{PS}
Consider dynamical percolation with parameter $p> p_c$ on the 
$d$-dimensional cubic lattice $\Z^d$. Then a.s\frenchspacing., for all times
$t$, there exists a unique infinite cluster.
\end{thm}

\begin{prop} \label{pr:treesupercr}\cite{PS}
Let $\Gamma$ be any infinite tree. If $p\in (p_c(\Gamma),1)$, 
then a.s\frenchspacing., 
there exist infinitely many infinite clusters for all $t$. 
\end{prop}

\subsection {Graphs with exceptional times exist.}

We now stick to the critical case. The following result says that
dynamical percolation can exhibit {\em exceptional times} which then
insures that this model is interesting.

\begin{thm} \label{th:G1G2}\cite{HPS}
There exists a $G_1$ which does not percolate at criticality
but such that there exist exceptional times of percolation. \\
There also exists a graph $\, G_2\, $ 
which percolates at criticality
but such that there exist exceptional times of nonpercolation.
\end{thm}

These graphs were obtained by replacing the edges of
$\Z^2$ by larger and larger finite 
graphs which simulated having the original
edge move back and forth quickly. These graphs have unbounded degree.
Many examples with bounded degree were later discovered and will be
discussed in later sections.

\subsection {High dimensional Euclidean case}

Once we have the existence of such examples, it is natural to ask what happens
on the standard graphs that we are familiar with.  The following result
answers this for the cubic lattice in high dimensions. (Recall
Hara and Slade \cite{HS} proved that ordinary percolation does not 
percolate at criticality in this case.)

\begin{thm}\cite{HPS} \label{th:Zd}
For the integer lattice $\Z^d$ with $d\ge 19,$ dynamical critical percolation
has no exceptional times of percolation.
\end{thm}

The key reason for this is a highly nontrivial result due to Hara and Slade, 
that says that if $\theta(p)$ is the probability that the origin percolates
when the parameter is $p$, then
\begin{equation} \label{eq:derive}
\theta(p)= O(p-p_c) \, .   
\end{equation}

In fact, Theorem \ref{th:Zd} was shown to hold under the assumption
(\ref{eq:derive}). 

\medskip\noindent
Outline of proof: 
We use a first moment argument together with the proof method of Proposition 
\ref{pr:noncrit}. 
We break the time interval $[0,1]$ into $n$ intervals each of 
length $1/n$. If we fix one of these intervals,
the set of edges which are open at {\em some time} during this interval
has density about $p_c+1/n$. Hence the probability that the
origin percolates with respect to these set of edges is by 
(\ref{eq:derive}) at most $O(1/n)$. It follows that the expected number of
these intervals where this occurs is at most $O(1)$. It can then be argued 
using Fatou's Lemma that a.s.\ there are at most finitely many
exceptional times during $[0,1]$ at which the origin percolates. To go from
there to no exceptional times can either be done by using some rather 
abstract Markov process theory or, as was done in the paper, by hand, which 
was not completely trivial. $\QED$

It is known, due to Kesten and Zhang \cite{KZ},
that (\ref{eq:derive}) fails for $\Z^2$. The question of whether
there are exceptional times for critical dynamical percolation on $\Z^2$ 
was left open in \cite{HPS}. (Recall there is no percolation at
a fixed time in this case.) This question will be resolved in Sections
\ref{sect:SS} and \ref{sect:GPS}.

\section {Exceptional times of percolation: the tree case}\label{sect:tree}

Although Theorem \ref{th:G1G2}
demonstrated the existence of graphs which have exceptional times
where the percolation picture is different from that at a fixed
time, in \cite{HPS}, a more detailed analysis was
done for trees. 

To explain this, we first need to give a little background from ordinary 
percolation on trees, results which are due to Lyons
and contained in \cite{Ly90} and 
\cite{Ly}. Lyons obtained an explicit condition
when a given tree percolates at a given value of $p$.
This formula becomes simplest when the (rooted)
tree is so-called
{\em spherically symmetric} which means that all vertices
at a given level have the same number of children, although
this number may depend on the given level. 

\begin{thm}\cite{Ly} \label{th:lyons}
(Case of spherically symmetric trees)\\
Let 0 be the root of our tree and 
$T_n$ be the number of vertices at the $n$th level. Then
the following hold. \\
(i) $p_c(T)=(\liminf_n T_n^{1/n})^{-1}$.\\
(ii) $\P_p(0\leftrightarrow \infty)>0$ if and only if
$\sum_n \frac{1}{w_n}<\infty$ where \\
$w_n:= \E_p[\mbox{number of vertices in $n$th level connected to 0}]
=p^n T_n$.
\end{thm}

Note that (i) essentially follows  from (ii).
In \cite{HPS}, an explicit condition was obtained 
for when, given a general tree and a value of $p$, there exists for
dynamical percolation a time at which percolation occurs. 
Again the spherically symmetric case is simplest.

\begin{thm}\cite{HPS} \label{th:HPStree}
(Case of spherically symmetric trees) \\
Let $w_n$ be as in the previous result. Then 
there exists a.s.\ a time at which percolation occurs
(equivalently with positive probability there exists a time in $[0,1]$
at which the root percolates) if and only if $\sum_n \frac{1}{nw_n}<\infty$.
\end{thm}

\medskip\noindent
Example: If $T_n \asymp 2^n n^{1/2}$, then the critical value is $1/2$,
there is no percolation at criticality but there are exceptional times
at which percolation occurs (which will have dimension $1/2$ by 
our next result). (The symbol $\asymp$ means that the ratio of the two 
sides is bounded away from zero and infinity by constants.)

\medskip\noindent
Outline of proof: For the ``if'' direction, we apply
the second moment method to the random variables
$$
Z_n:=\sum_{v \in T_n} U_v
$$
where $U_v$ is the Lebesgue amount of time during $[0,1]$ 
that the root is connected to $v$ and
$T_n$ now denotes the set of vertices at the $n$th level. 
(The ``second moment method'' means that one computes both the first
and second moments of a given nonnegative random variable $X$
and applies the Cauchy-Schwarz inequality to obtain
$P(X >0) \ge (E[X]^2)/E[X^2]$.)
For the ``only if'' direction, we perform
a slightly nontrivial {\em first space-time decomposition argument} 
as follows.
If there is some $v$ in $T_n$ which percolates to the root 
at some time $t$ during $[0,1]$, we show that one can pick such a random pair 
$(v',t')$ in such a way that $E[Z_n\mid (v',t')]$ is much larger than
$E[Z_n]$. This implies that it is very unlikely that such a pair exists.
$\QED$

In \cite{HP}, it was shown that some of the more refined results
in \cite{Ly} and Theorem \ref{th:HPStree} above together yield that
a spherically symmetric tree has an exceptional time of
percolation if and only if $\int_{p_c}^1 1/\theta(p) < \infty$.
(As we will see, this equivalence also turns out to be true for 
$\Z^d$ with $d=2$ or $d\ge 19$.)
It was also shown in \cite{HP} that for general trees,
the ``only if'' holds but not the ``if''.

In the spherically symmetric case, the Hausdorff dimension was also
determined in \cite{HPS}.

\begin{thm}\cite{HPS} \label{th:HPStreeHD}
Consider a spherically symmetric tree and 
let $w_n$ be as in the two previous results. Then the 
Hausdorff dimension of the set of times at which percolation occurs is given by
\[
\sup \Big\{ \, \alpha \in [0,1] \: :  \: \sum_{n=1}^\infty 
  \frac{n^{\alpha-1}}{w_n} \, < \infty \, \Big\} \, .
\]
\end{thm}

While we do not state the results here, we mention that
\cite{PS} went a good deal further for 
spherically symmetric trees. What was obtained
in \cite{PS} was (see Corollary 5.1 there) a ``capacity 
condition'' for which subsets of
time have the property that with positive probability
they contain a percolating time.  This is analogous
to the well-known Kakutani criterion (in terms of
Newtonian capacity) of which subsets in Euclidean
space intersect a Brownian motion path with 
positive probability. 

Once we have such a capacity condition, Peres'
intersection equivalence theory (\cite{Pe} and \cite{Pe2})
leads to a criterion for when there are exceptional times at which
there are at least $k$ infinite clusters. We can 
in particular construct trees for which there are
times at which we have (say) 6 infinite clusters but no times 
at which there are 7 infinite clusters. In addition,
various Hausdorff dimensions of these different exceptional time
sets can be computed.

On a much less formal note, my personal feeling is that the set of 
exceptional times in these cases, very vaguely speaking, 
might have a similar structure to the set of so-called ``slow'' 
points for Brownian motion.
In Section~\ref{sect:PSS}, we will see a very different type of
set of exceptional times for dynamical percolation
and I believe that in this latter
case, this set might behave more like the 
set of so-called ``fast'' points for Brownian motion. Here are two words
explaining this vague connection (for those who know these concepts
from Brownian motion). 
For a time point $s$ to be a slow point,
it must be the case that $|B(t+s)-B(s)|$ does not go above a certain
(well-specified) level for {\em all} values of small $t$ while
a time $s$ is an exceptional time of percolation if the origin percolates out
to {\em all} distances. 
On the other hand, for a time point $s$ to be a fast point,
it must be the case that $|B(t+s)-B(s)|$ goes above a certain
(well-specified) level for an {\em infinite number} of arbitrarily
small $t$ while
a time $s$ is an exceptional time of non-percolation say for a tree
if we have an {\em infinite number} of cut-sets which are off. This type
of structure looks like something which is called a limsup fractal;
see \cite{KPX}.

We end by mentioning that in \cite{K}, Khoshnevisan extended 
to general trees the result of Peres and Steif 
determining which time sets contain percolating times.
Results concerning the Hausdorff dimension of exceptional times
for general trees are also obtained in \cite{K}. Some of the techniques 
use methods from \cite{Kbook}.

\section {Noise sensitivity, noise stability,
influence and the Fourier spectrum}\label{sect:NS}

The study of noise sensitivity and noise stability 
for Boolean functions was initiated in \cite{BKS}. 
The concepts discussed in this section come from or are motivated by
this source. See \cite{Othesis} for a nice survey of 
noise sensitivity and its applications in
theoretic computer science. See also \cite{KS} for related matters. 

\subsection {Definition of noise sensitivity, noise stability 
and some examples} \label{subsect:4.1}

As indicated
in the introduction, the notion of noise sensitivity is {\em both}
an interesting concept in itself {\em and} what is needed to carry out
the second moment arguments necessary for the results
described in Sections \ref{sect:SS} and \ref{sect:GPS}.

Let $\omega$ be uniformly chosen from
the $n$-dimensional discrete cube $\{0,1\}^n$ (which we can think of
as $n$ fair coin flips) and let 
$\omega_\epsilon$ be $\omega$ but with each bit independently
``rerandomized'' with probability $\epsilon$. 
``Rerandomized'' means the bit (independently of everything
else) rechooses whether it is 1 or 0, each with probability $1/2$.
(Of course $\omega_\epsilon$ has the same distribution as $\omega$).

Let $f:\{0,1\}^n \rightarrow \{\pm 1\}$ or $\{0,1\}$ be arbitrary.
We are interested in the covariance between
$f(\omega)$ and $f(\omega_\epsilon)$. In most cases of interest,
we will have a sequence $\{f_n\}$ where 
$f_n:\{0,1\}^{m_n} \rightarrow \{\pm 1\}$ or $\{0,1\}$
and we are interested in the asymptotic behavior of the covariance above.

\begin{defn}
The sequence $\{f_n\}$ is {\bf noise sensitive} if for every $\epsilon> 0$,
$$
\lim_{n\to\infty} \E[f_n(\omega)f_n(\omega_\epsilon)]- \E[f_n(\omega)]^2 =0.
$$
\end{defn}

Usually $f$ is an indicator of an event $A$
and this then says that the
two events $\{\omega\in A\}$ and $\{\omega_\epsilon\in A\}$
are close to independent for $\epsilon$ fixed and $n$ large.
The following notion captures the opposite situation where
the two events above are close to being the same event 
if $\epsilon$ is small, uniformly in $n$.

\begin{defn}
The sequence $\{f_n\}$ is {\bf noise stable} if 
$$
\lim_{\epsilon\to 0} 
\sup_n \P(f_n(\omega)\neq f_n(\omega_\epsilon))=0.
$$
\end{defn}

It is easy to check that $\{f_n\}$ is both noise sensitive and noise stable 
if and only if the sequence of variances $\{\Var(f_n)\}$ goes to 0.
Here are two easy examples where $m_n$ is taken to be $n$.

\medskip\noindent
Example 1:
$f_n(\omega)=\omega_1$ (i.e., the first bit).

\medskip\noindent
Example 2:
$f_n(\omega)$ is the parity of the number of 1's in $\omega$.

\medskip\noindent
It is easy to check that Example 1 is noise stable while Example 2 is 
noise sensitive. We will see later why these two examples are the two extreme
examples. A more interesting example is the following which, while
it is not immediately obvious, turns out to be noise stable
as shown in \cite{BKS}.

\medskip\noindent
Example 3 (Majority): Let $m_n=2n+1$. Let
$f_n(\omega)$ be 1 if there is a majority of 1's and
0 if there is a majority of 0's.

\medskip\noindent
A much more interesting example is the following. 

\medskip\noindent
Example 4: Let $f_n$ be the indicator
function of a left to right crossing of the box $[0,n]\times[0,n]$
for critical percolation either for the ordinary lattice
$\Z^2$ or for the hexagonal lattice. For the hexagonal lattice,
the box has to be slightly modified so that it is a union of hexagons.
(Note $m_n$ is of order $n^2$ in this case.) 

\begin{thm} \label{th:BKS}\cite{BKS}
The sequence $\{f_n\}$ in Example 4 is noise sensitive.
\end{thm}

\noindent
In fact, it was shown that
\begin{equation}\label{e.z19}
\lim_{n\to\infty}\E[f_n(\omega)f_n(\omega_{\epsilon_n})]- \E[f_n(\omega)]^2=0
\end{equation}
even when $\epsilon_n$ goes to 0 with $n$ provided that
$\epsilon_n\ge C/\log(n)$ for a sufficiently large $C$. Clearly for any
sequence of Boolean functions, if $\epsilon_n$ goes to 0 sufficiently fast,
we have
\begin{equation}\label{e.stable}
\lim_{n\to \infty} 
\P(f_n(\omega)\neq f_n(\omega_{\epsilon_n}))=0
\end{equation}
for the trivial reason that in that case
$\P(\omega\neq \omega_{\epsilon_n})\to 0$.

\subsection {The noise sensitivity exponent, the noise stability exponent,
and influence.}

For sequences $\{f_n\}$ which are noise sensitive, it might be hoped
that (\ref{e.z19}) is still true when $\epsilon_n$ decays as some
power of $1/n$ and this was explicitly asked in
\cite{BKS} for crossings in critical percolation.
This suggested ``critical exponent'' is what we call the noise 
sensitivity exponent. The following definitions now seem natural.

\begin{defn}
The {\bf noise sensitivity exponent} $(\SENS(\{f_n\}))$ 
of a sequence $\{f_n\}$ is defined 
to be
$$
\sup\{\alpha: \mbox{ (\ref{e.z19}) holds with } 
\epsilon_n=(1/n)^\alpha\}.
$$ 
The {\bf noise stability exponent} $({\rm STAB}(\{f_n\}))$ 
of a sequence $\{f_n\}$ is defined to be
$$
\inf\{\alpha: \mbox{ (\ref{e.stable}) holds with } 
\epsilon_n=(1/n)^\alpha\}.
$$
\end{defn}

\medskip\noindent {\em Remarks:} \\
\noindent
1. We will see later that 
$\E[f(\omega)f(\omega_\epsilon)]- \E[f(\omega)]^2$ is nonnegative and 
decreasing in $\epsilon$. It easily follows that
$\P(f_n(\omega)\neq f_n(\omega_{\epsilon}))$ is increasing in $\epsilon$. 
From this, it is easy to see that
$\SENS(\{f_n\})\le \STAB(\{f_n\})$ unless the
variances $\Var(\{f_n\})$ go to 0. \\
\noindent
2. One might hope that $\SENS(\{f_n\})= \STAB(\{f_n\})$. 
First, this can fail for trivial reasons
such as the $f_n$'s for even $n$ might have nothing to do 
with the $f_n$'s for odd $n$. However, this may fail
for more interesting reasons. 
Using $\STAB(\{A_n\})$ for $\STAB(\{I_{A_n}\})$ and similarly for
$\SENS$, if, for example, $A_n$ and $B_n$ are
independent for each $n$, have probabilities near $1/2$ and satisfy
$$
\SENS(\{A_n\})= \STAB(\{A_n\})=a<b=\SENS(\{B_n\})= \STAB(\{B_n\}),
$$
then it is easy to check that
$$
\SENS(\{A_n\cap B_n\})=a < b=\STAB(\{A_n\cap B_n\}).
$$
In such a case, for $\epsilon_n=(1/n)^c$ with $c\in (a,b)$,
the correlation between $\{\omega\in A_n\cap B_n\}$ 
and $\{\omega_{\epsilon_n}\in A_n\cap B_n\}$ neither goes to 0 nor to being
perfectly correlated. The question of under what conditions
we have $\SENS(\{f_n\})= \STAB(\{f_n\})$ turns out to be a fairly subtle 
one. 
\\
\noindent
3.  If $m_n\asymp n^\sigma$, then it is trivial to check that 
$\STAB(\{f_n\})  \le \sigma$ since if
$\epsilon_n=(1/n)^{\sigma+\delta}$ for some fixed $\delta$, then
$\P(\omega\neq \omega_{\epsilon_n})\to 0$.
\\
\noindent
4. It is natural to ask for bounds on these exponents for general Boolean 
functions; this will be discussed at the end of subsection
\ref{subsection:Fourier}.

\medskip\noindent
The next important notion of total influence will give us a more interesting 
upper bound on the noise stability exponent than that provided by comment 3 
above.

\begin{defn}
Given a function  $f:\{0,1\}^n \rightarrow \{\pm 1\}$ or $\{0,1\}$,
let $I_i(f)$, which we call the {\bf influence} of the $i$th variable on $f$, 
be the probability that all the variables other than the
$i$th variable do not determine $f$. I.e., letting $\omega^i$ be
$\omega$ but flipped in the $i$th coordinate,
$$
I_i(f):=\P(\omega: f(\omega)\neq f(\omega^i)).
$$
The {\bf total influence}, denote by $I(f)$, is defined to be $\sum_i I_i(f)$.
If $f(\omega)\neq f(\omega^i)$, we say that $i$ is {\bf pivotal} for $f$
and $\omega$ and hence  $I(f)$ is the expected number of pivotal bits of $f$.
\end{defn}

\medskip\noindent
The following I believe was ``known'' in the community. The argument
below however I first saw given by Christophe Garban in the context of 
percolation.

\begin{thm} \label{th:influence}
Consider a sequence $f_n:\{0,1\}^{m_n} \rightarrow \{\pm 1\}$ or $\{0,1\}$
and assume
that $I(f_n)=n^{\rho+o(1)}$. Then $\STAB(\{f_n\})  \le \rho$.
\end{thm}

\medskip\noindent
Proof:
We need to show that if $\alpha > \rho$ and 
$\epsilon_n=(1/n)^\alpha$, then (\ref{e.stable}) holds.
Let $\omega_0,\omega_1,\ldots,\omega_{k_n}$ be such that one obtains
$\omega_{i+1}$ from $\omega_i$ by choosing independently a bit at random
and rerandomizing it. It is immediate that
$\P(f_n(\omega_i)\neq f_n(\omega_{i+1}))=\frac{I(f_n)}{2m_n}$ 
from which one gets
$\P(f_n(\omega_0)\neq f_n(\omega_{k_n}))\le \frac{k_nn^{\rho+o(1)}}{m_n}$.
If $\epsilon_n=(1/n)^\alpha$, then $\omega_{\epsilon_n}$ is sort of like
$\omega_{k_n}$ with $k_n=\frac{m_n}{n^\alpha}$ and so
$$
\P(f_n(\omega)\neq f_n(\omega_{\epsilon_n}))\sim
\P(f_n(\omega_0)\neq f_n(\omega_{k_n}))\le
\frac{m_nn^{\rho+o(1)}}{m_nn^\alpha}
$$
which goes to 0 as $n\to\infty$ if $\alpha > \rho$.
The ``sort of like'' above and the imprecise $\sim$ are trivial to make
rigorous and correct using standard large deviations for binomial random 
variables. 
$\QED$

\medskip\noindent {\em Remarks:} \\
An example where we have strict inequality
is the ``majority function'' of Example 3. It is easy to check
that for this example $\rho$ in Theorem \ref{th:influence} is $1/2$ but,
by the noise stability of this example, we have that
$\STAB(\{f_n\})=0$. 
One explanation of the failure of having a converse to Theorem 
\ref{th:influence} in this case is that the expected number of pivotals 
is not so relevant:
the random number  $N_n$ of pivotals is not at all concentrated around its
mean $\E[N_n]=I(f_n)\asymp n^{1/2}$ but rather it goes to 0 in probability.

We will see an alternative proof of Theorem \ref{th:influence} 
using Fourier analysis in the next subsection.

Of the very many nice results in \cite{BKS}, we mention the following
one which almost characterizes noise sensitivity in terms of influences.
As the authors mention, this result could be used to prove Theorem
\ref{th:BKS} but they instead use a different approach.
A function $f$ is {\bf montone} if $x\le y$ (meaning
$x_i\le y_i$ for each $i$) implies that $f(x)\le f(y)$. The proof of the
following result uses the Fourier spectrum (see the next subsection) and
a notion known as hypercontractivity.

\begin{thm} \label{th:ns=inf}\cite{BKS}
Consider a sequence of Boolean functions $\{f_n\}$.
If \newline
$\lim_{n\to\infty}\sum_i I_i(f_n)^2=0$, then the sequence
$\{f_n\}$ is noise sensitive. The converse is true if the $f_n$'s
are monotone. (Example 2 shows that monotonicity is needed for the converse.)
\end{thm}

We end this section by going back to Example 4. It was asked
in \cite{BKS} whether the noise sensitivity
exponent for this sequence is 
$3/4$. The heuristic for this guess is the following. An edge 
(or hexagon) is pivotal if and only if there are 4 disjoint 
monochromatic paths
alternating in color from the edge (or hexagon) to the top, right, bottom
and left sides of the box. This is close to the 4-arm exponent
which by \cite{SW} was proved to be behave like $n^{-5/4}$. If boundary terms
do not matter much, $I(f_n)$ should then be about $n^{3/4}$.
One might hope that the number $N_n$ of pivotals is quite concentrated
around its mean $I(f_n)$; in this direction, it is well-known for example
that $\E[N_n^2]=O(1)\E[N_n]^2$. This gives hope that Theorem \ref{th:influence} 
is now tight in this case and that $3/4$ might be both the noise 
sensitivity and noise stability exponents. 
In Section \ref{sect:GPS}, we will see that this is indeed the case,
but the proof will not go through understanding pivotals but
rather through Fourier analysis. This brings us to our next topic.

\subsection {The Fourier spectrum} \label{subsection:Fourier}

It turns out that the best and proper way to analyze the above problems is to 
use Fourier analysis. The set of all functions $f:\{0,1\}^n \rightarrow \R$ 
is a $2^n$ dimensional vector space with orthogonal basis 
$\{\chi_S\}_{S\subseteq \{1,\ldots,n\}}$ where
$$
\chi_S(\omega^1,\ldots,\omega^n):=\left\{ \begin{array}{ll}
-1\quad&\text{if \# of 1's in
$\{\omega^i\}_{i\in S}$ is odd}
\\1\quad&\text{if \# of 1's in
$\{\omega^i\}_{i\in S}$ is even. }
\end{array} 
\right.
$$
So $\chi_\emptyset \equiv 1$. We then can write
$$
f:=\sum_{S\subseteq \{1,\ldots,n\}}  \hat{f}(S) \chi_S.
$$
(In fancy terms, the various $\chi_S$'s
are the so-called characters on the group $\Z_2^n$
but everything below will be from first principles).
The $\hat{f}(S)$'s are called the Fourier coefficients.

The reason that $\{\chi_S\}$ is a useful basis is that
they are eigenfunctions for the discrete time Markov chain
which takes $\omega$ to $\omega_\epsilon$.
It is an easy exercise to check that
$$
\E[\chi_S(\omega_\epsilon)|\omega]=
(1-\epsilon)^{|S|}\chi_S(\omega)
$$
and that
\begin{equation} \label{eq:expansion}
\E[f(\omega)f(\omega_\epsilon)]= 
\E[f(\omega)]^2+ \sum_{k=1}^n  (1-\epsilon)^k
\sum_{|S|=k}\hat f(S)^2.
\end{equation}

This formula which first (I believe) appeared in \cite{BKS}
makes clear the central role played
by the Fourier coefficients with respect to
questions involving noise. Note that we see the nonnegativity  
of the covariance between $f(\omega)$ and $f(\omega_\epsilon)$
and that it is decreasing in $\epsilon$ as we had claimed 
earlier. (I am not aware of any coupling proof of this latter fact;
we can of course couple $\omega$,
$\omega_\epsilon$ and $\omega_{\epsilon+\delta}$ so that 
$\omega_\epsilon$ agrees with $\omega$ in more places than
$\omega_{\epsilon+\delta}$ does but in view of Example 2 
in Subsection \ref{subsect:4.1}, this
does not seem to help.) Crucially, we see that the 
covariance between $f(\omega)$ and $f(\omega_\epsilon)$ 
is small when most of the ``weight'' of 
these coefficients are on the $S$'s with larger cardinality while
the covariance is largest when most of the 
``weight'' of these coefficients are on the smaller $S$'s.

In Subsection \ref{subsect:4.1},
Example 1 is the function $(1-\chi_{\{1\}})/2$
while Example 2 is the function $\chi_{\{1,\ldots,n\}}$
from which we now see why these are extreme examples 
as we mentioned earlier.

We now restrict to $f$'s which take values in $\{-1,1\}$.
The advantage in doing this is that we have (due to the
Pythagorean theorem or Parseval's formula)
$$
\sum_{S\subseteq \{1,\ldots,n\}}\hat f(S)^2=1.
$$
(For those who do not like the restriction of $\{-1,1\}$
on the range of $f$ since you are interested in indicator
functions of events, you can just consider the function which is 1 on
the event in question and $-1$ on its complement and then easily translate 
back the results below in terms of your events.)

Given the above, we can let 
${\mathcal S}$ be a random subset of $\{1,\ldots,n\}$
given by $\P({\mathcal S}=S)=\hat f(S)^2$. The idea of looking at this as a 
probability distribution on the subsets of $\{1,\ldots,n\}$
was proposed in \cite{BKS} and is called the 
{\bf spectral sample} or {\bf spectral measure} of $f$.  (We will
not be careful to distinguish between ${\mathcal S}$ and its distribution.)
(\ref{eq:expansion}) can now be written as
(with the two expectations being on different spaces)
\begin{equation} \label{eq:expansionrrandom}
\E[f(\omega)f(\omega_\epsilon)]= \E[(1-\epsilon)^{|\mathcal S|}].
\end{equation}
This equation demonstrates the important fact that a sequence 
$\{f_n\}$ is noise sensitive if and only if the corresponding
spectral measures $\{\mathcal S_n\}$ satisfy
$|\mathcal S_n|\to \infty$ in distribution provided we remove the point mass
at $\emptyset$ (which corresponds to subtracting the squared mean).

We now give an alternative proof of Theorem \ref{th:influence} which
is taken from \cite{GPS}.

\newpage

\medskip\noindent
Proof \cite{GPS}:
We need only one thing which we do not prove here which 
comes from \cite{KKL};
one can also see a proof of this in \cite{GPS}. This is that 
\begin{equation} \label{eq:inf=S}
\E[|\mathcal S|]=I(f).
\end{equation}
It clearly suffices to show that 
if $\alpha > \rho$ and $\epsilon_n=(1/n)^\alpha$, then
$$
\lim_{n\to\infty}\E[f_n(\omega)f_n(\omega_{\epsilon_n})]=1.
$$
(\ref{eq:expansionrrandom}), Jensen's inequality and
$\E[|{\mathcal S}_n|]=I(f_n)$ yields
$$
\E[f_n(\omega)f_n(\omega_{\epsilon_n})]=
\E[(1-\epsilon_n)^{|{\mathcal S}_n|}]\ge
(1-\epsilon_n)^{\E[|{\mathcal S}_n|]}=
(1-\epsilon_n)^{I(f_n)}
$$
$$
=
(1-(1/n)^\alpha)^{n^{\rho+o(1)}}.
$$
This goes to 1 since $\alpha > \rho$.
$\QED$

In the context of Theorem \ref{th:influence}, it would be nice to know under
which further conditions we could conclude the reverse inequality and
even that $\SENS(\{f_n\})  \ge \rho$. One sufficient condition is that 
the distributions of $|{\mathcal S}_n|$, normalized by their means and with the
point mass at 0 removed is tight on $(0,\infty)$. 
(Tightness at $\infty$ follows from Markov's inequality;
the key point is tightness near 0.)
We make this precise in the following result. This result is proved
by just abstracting an easy small part of an argument from \cite{GPS}. 

\begin{thm} \label{thm:meta} 
Assume $f_n:\{0,1\}^{m_n} \rightarrow \{\pm 1\}$ and that $I(f_n)=n^{\rho+o(1)}$.
Let  ${\mathcal S}_n$ be the spectral sample corresponding to $f_n$.
Assume that for every $\gamma>0$, there is $\delta>0$ so that for all $n$
\begin{equation}\label{eq:tightness}
\P(|{\mathcal S}_n| < \delta \E[|{\mathcal S}_n|],{\mathcal S}_n\neq \emptyset)
<\gamma.
\end{equation} 
Then $\SENS(\{f_n\}) \ge \rho$.
\end{thm}

\medskip\noindent
Proof:
We need to show that if $\alpha < \rho$, 
(\ref{e.z19}) holds when $\epsilon_n=(1/n)^\alpha$.
The difference in (\ref{e.z19}) is by (\ref{eq:expansion}) and
(\ref{eq:expansionrrandom})
simply $\E[(1-\epsilon_n)^{|\mathcal S_n|}I_{{\mathcal S}_n\neq \emptyset}]$.
Fix $\gamma >0$ and choose $\delta$ as in the assumption. We have
$$
\E[(1-\epsilon_n)^{|\mathcal S_n|}I_{{\mathcal S}_n\neq \emptyset}]
\le \P(|{\mathcal S}_n| < \delta \E[|{\mathcal S}_n|],{\mathcal S}_n\neq 
\emptyset) + (1-\epsilon_n)^{\delta \E[|\mathcal S_n|]}.
$$
The first term is at most $\gamma$ for all $n$ by the choice of $\delta$
and the last term goes to 0 since $\delta$ is fixed,
$\alpha < \rho$ and using (\ref{eq:inf=S}).  
Since $\gamma$ was arbitrary, we are done.
$\QED$

It is now interesting to look again at the majority function
for which $\rho=1/2$ but $\SENS(\{f_n\}) =\STAB(\{f_n\})=0$. In this case,
the spectral measures do not satisfy the necessary
tightness condition above but rather these distributions
normalized by their means approach the point mass at 0.
This follows from the known noise stability of majority,
the fact (see Theorem 1.9 in \cite{BKS}) that stability in general implies
tightness at $\infty$ for the {\em unnormalized} spectral measures
and the fact that the expected value of the spectral size is going to
$\infty$. See \cite{Othesis} for details concerning the spectral 
measure of majority.

The following is an exercise for the reader which relates the lower tail
of the spectral measures with noise sensitivity;
Theorem \ref{thm:meta} already gave some relationship between these.

\medskip\noindent
Exercise 1: Let $f_n$ be an arbitrary sequence of Boolean functions
with corresponding spectral samples ${\mathcal S}_n$. \\
(i) Show that $\P(0<|{\mathcal S}_n| \le A_n)\rightarrow 0$ implies that
$\E[(1-\epsilon_n)^{|\mathcal S_n|}I_{{\mathcal S}_n\neq \emptyset}]\rightarrow 0$ 
if $\epsilon_n A_n \rightarrow \infty$.  \\
(ii) Show that
$\E[(1-\epsilon_n)^{|\mathcal S_n|}I_{{\mathcal S}_n\neq \emptyset}]\rightarrow 0$ 
implies that
$\P(0<|{\mathcal S}_n| \le A_n)\rightarrow 0$ if $\epsilon_n A_n =O(1)$.\\
In particular, $\SENS(\{f_n\}) \ge \alpha$ if and only if
$\P(0<|{\mathcal S}_n| \le n^{\alpha-\delta})\rightarrow 0$ for all
$\delta > 0$.

\medskip\noindent
Exercise 2: Show that $\STAB(\{f_n\}) \le \alpha$ if and only if
$\P(|{\mathcal S}_n| \ge n^{\alpha+\delta})\rightarrow 0$ for all
$\delta > 0$.

We end this subsection with a brief discussion of general upper bounds on
our noise sensitivity and noise stability exponents. We stick for simplicity
and without loss of generality to $m_n=n$. We have seen then that 1 is an
upper bound for these exponents. On the other hand, the parity function,
Example 2, is easily seen to have both these exponents being 1. 
This question turns out to be much more interesting if we restrict
to the important subclass of {\bf monotone} Boolean functions.

\begin{thm} \label{thm:monotone} 
Assume $f_n:\{0,1\}^{n} \rightarrow \{\pm 1\}$ be monotone. 
Then $\STAB(\{f_n\}) \le 1/2$.
\end{thm}

\medskip\noindent
Outline of proof: It is an exercise to check, crucially using
the monotonicity, that $I_i(f_n)=  |\hat f_n(\{i\})|$. It follows that
$\sum_i I_i(f_n)^2\le 1$. The Cauchy-Schwartz inequality now yields that
$\sum_i I_i(f_n)\le \sqrt{n}$. The result now follows from
Theorem \ref{th:influence}.
$\QED$

\medskip\noindent {\em Remark:} \\
The above theorem can be proved (and has been) without the use
of Fourier analysis and is known in various contexts.

\medskip\noindent
Answering a question in \cite{BKS}, it was 
shown in \cite{MO} that the above result is optimal by giving a
sequence $\{f_n\}$ of monotone functions with $\STAB(\{f_n\})=1/2$. 
By tweaking these examples, one can also obtain a sequence with 
$\SENS(\{f_n\})=1/2$.

\section {Critical exponents for percolation} 
\label{sect:crit}
The exact values for critical exponents for percolation
on the hexagonal lattice are a crucial ingredient in
the study of exceptional times for dynamical percolation
and noise sensitivity for crossing probabilities. We therefore
briefly describe these.

Let $A^k_R$ be the event (for ordinary percolation)
that there are $k$ disjoint 
monochromatic paths from within distance (say) $2k$ of the 
origin all of which reach distance $R$ from the origin and such that they
are not all of the same color. This is referred to as the $k$-arm event.
Let $A^{k,H}_R$ be the analogous event but where we restrict to the upper half 
plane and where the restriction ``not all of the same color'' may be dropped.
This is referred to as the half-plane $k$-arm event. All these events
decay as powers of $R$ and the exact power is called the corresponding
critical exponent.

\begin{thm}\label{th:critexp} For the hexagonal lattice, we have \\
(i) \cite{LSW} $\P(A^1_R)= R^{-5/48+o(1)}$ \\
(ii) \cite{SW} For $k\ge 2$, $\P(A^k_R)= R^{-(k^2-1)/12+o(1)}$ \\
(iii) \cite{SW} For $k\ge 1$, $\P(A^{k,H}_R)= R^{-k(k+1)/6+o(1)}$ \\
(iv)  \cite{SW} $\theta(p)=(p-1/2)^{5/36+o(1)}$
\end{thm}

It was shown by Kesten \cite{Kscaling} that (iv) follows from (i) and the 
case $k=4$ in (ii). See \cite{Nolin} for a detailed derivation of this.

For dynamical percolation on the hexagonal lattice,
$A^1_R$, $A^2_R$ and $A^{1,H}_R$ and their exact critical exponents
were relevant in the work in \cite{SS} while
$A^1_R$, $A^4_R$, $A^{1,H}_R$ and even some ``corner-plane'' events
and their exact critical exponents 
were relevant in the work in \cite{GPS}.

We finally mention that the above proofs rely on the conformal invariance
of percolation on the hexagonal lattice proved by Smirnov (see 
\cite{Smirnov}) and the convergence of the discrete interface in
critical percolation to $SLE_6$ (see \cite{CN} and \cite{Smirnov}).
SLE originally stood for stochastic L\"{o}wner evolution when it 
was introduced by Schramm in \cite{Schramm} and is presently called the
Schramm-L\"{o}wner evolution. It has one parameter, usually called
$\kappa$, and therefore written SLE$_{\kappa}$. As one varies
$\kappa$, these yield random curves which describe many 2-dimensional
critical systems.

\section {Exceptional times and positivity of the noise 
sensitivity exponent for the hexagonal lattice} \label{sect:SS}

Considering the question of whether there are exceptional times
for dynamical percolation in $\Z^2$, while this was not accomplished 
in \cite{SS}, it was proved in this paper that exceptional times do exist
for the hexagonal lattice. What allowed the
proof to go through for the hexagonal lattice is that
various exact critical exponents from critical percolation 
have been established for the hexagonal lattice. These same critical
exponents are expected to hold for $\Z^2$ but have not at this point
been established.

\begin{thm} \label{th:SS}\cite{SS} 
For critical dynamical percolation on the hexagonal lattice, there exist 
exceptional times of percolation and the Hausdorff dimension of 
the set of such times is in $[1/6,31/36]$.
\end{thm}

As far as the noise sensitivity exponent for left to right crossing of an 
$n\times n$ square, the following was shown.

\begin{thm} \label{th:SSNoise}\cite{SS} 
Consider the sequence $\{f_n\}$ of indicator functions 
for a left to right crossing of an $n\times n$ square in the 
hexagonal lattice. Then $\SENS(\{f_n\}) \ge 1/8$.
For the square lattice, $\SENS(\{f_n\}) > 0$.
\end{thm}

This was the first result where one obtained a positive
lower bound on $\SENS(\{f_n\})$ for crossing probabilities. 
One of the key steps in \cite{SS} is the following result which
gives conditions under which one can obtain some bounds on the Fourier
coefficients. We will not give any indication of its proof here. We hope
that this result will be applicable in other contexts in order to bound
the ``level-$k$'' Fourier coefficients. 

\begin{thm} \label{th:algorithm} \cite{SS} \\
Let $f:\{0,1\}^n\to \R$. Let $A$ be a randomized algorithm 
determining the value of $f$. This means that $A$ examines the 
input bits of $f$ one by one, where the choice of the next bit 
examined may depend on the values of the bits examined so far and
on some additional exterior randomness. The algorithm $A$ may of course stop 
once it knows the output. Let $J\subseteq \{1,2,\dots,n\}$
be the (random) set of bits examined by the algorithm. 
(Crucially, $J$ need not be all the bits since based on the bits
examined at some point, the output might at that point be determined.)
Set 
$$
\delta_A:=\sup \{ \P(i\in J):i\in\{1,2,\dots,n\}\}.
$$
Then, for every $k=1,2,\dots$,
the Fourier coefficients of $f$ satisfy
$$
\sum_{|S|=k}\hat f(S)^2  \le \delta_A \,k\,\|f\|^2,
$$
where $\|f\|$ denotes the $L^2$ norm of $f$.
\end{thm}

\medskip\noindent
We first give an outline of the proof of Theorem \ref{th:SSNoise}.

\medskip\noindent
Outline of proof: \\
Step 1: Find a randomized algorithm to detect whether 
there is a left to right crossing of whites such that a fixed 
hexagon is looked at with probability at most $(1/n)^{1/4+o(1)}$. 
If we place white hexagons on the right and top sides of the box and
black hexagons on the bottom and left sides of the box, we can start at
the bottom right and follow the path which always keeps white on the
right side and blacks on the left. This is called an {\bf interface}.
This interface will end up either hitting the left side before the top,
which implies there is a left to right white crossing or it will 
end up either hitting the top side before the left,
which implies there is no left to right white crossing.
If we ``follow'' this interface, revealing hexagon colors as we need to
know how the interface develops, this will yield a randomized
algorithm which will determine if there is a crossing.
In addition, hexagons near the center will be 
looked at with probability at most $O(1)(1/n)^{1/4+o(1)}$
since to be looked at, one must see the
``2-arm'' event emanating from that hexagon and one can apply 
Theorem \ref{th:critexp}(ii). This does not work however
for hexagons near the boundary. To get an 
algorithm which looks at {\em every} hexagon with the 
above probability, one does some random modification
of the above where one runs two interfaces from a random point on
the right side. The argument then requires using the 1-arm 
half plane exponent as well. \\
Step 2:  Step 1 and Theorem \ref{th:algorithm} gives us a bound
on the sum of the ``level-$k$'' Fourier coefficients. We plug that
into (\ref{eq:expansion}) and compute. The $\Z^2$ case is similar but 
there we do not have the explicit critical exponents at our disposal.
$\QED$

\medskip\noindent {\em Remarks:} \\
We might hope that one can bring the value $1/8$
up to the suggested value of $3/4$
by finding better algorithms to which
we can apply Theorem \ref{th:algorithm}.
However, this is not possible. As mentioned in \cite{SS},
a general inequality of O'Donnell and Servedio 
or Theorem \ref{th:algorithm} applied in the case $k=1$ allows us to conclude 
that any algorithm will have a $\delta$ of at least  
$(1/n)^{1/2+o(1)}$. The existence of such an algorithm would (as above)
bring the value of $1/8$ up to $1/4$ and hence the best this method 
could yield is a noise sensitivity exponent of $1/4$ 
(unless one improves Theorem \ref{th:algorithm} itself).
It is worth pointing out here that an algorithm which is conjecturally
better than the one given in the proof of Theorem \ref{th:SSNoise}
is the one where the hexagon chosen at a given
time is the one with the largest influence at the time. This is related
to playing random-turn hex; see \cite{PSSW}.

\medskip\noindent
We now give an outline of the proof of Theorem \ref{th:SS}.

\medskip\noindent
Outline of proof: \\
We first explain the existence of exceptional times.
We let $X_R:=\int_0^1 1_{V_{t,R}} \, dt$ where
$V_{t,R}$ is the event that at time $t$ there is an open path 
from the origin to distance $R$ away. The key step is to show that
$\E[X_R^2]\le O(1)\E[X_R]^2$ (i.e., we use the second moment method)
since from here the result is standard. Note the first moment is
just $\P(A^1_R)$ from Section \ref{sect:crit}.
The key to getting a good bound on the second moment 
is to get a good bound on $\P(V_{t,R}\cap V_{0,R})$. 
Using independence, we see this is at most
$\P(V_{0,r})\P(V_{t,r,R}\cap V_{0,r,R})$ where
$V_{t,r,R}$ is the event that at time $t$ there is an open path from 
distance $r$ away to distance $R$ away. We note that looking
at our process at times $0$ and $t$ is exactly looking at a configuration
and its noisy version, which we studied earlier, with
$\epsilon$ being $1-e^{-t}$. For the second factor,
we use (\ref{eq:expansion}), construct a randomized algorithm for 
this event with a good $\delta$ (which is somewhat harder than
in Theorem \ref{th:SSNoise}) and apply Theorem \ref{th:algorithm} 
to bound the relevant Fourier coefficients. The rest is algebra.

For the Hausdorff dimension, the lower bound is obtained using the
calculation in the first paragraph together with Frostman's Theorem. 
The upper bound is easier; we use the method of proof of 
Theorem \ref{th:Zd} together with Theorem \ref{th:critexp}(iv).
$\QED$

In \cite{SS}, other types of events are also looked at such as $k$-arm
events in wedges and cones and upper and lower bounds on the Hausdorf
dimension of sets of exceptional times are obtained. 
These upper and lower bounds are however never matching.

\section {The exact Hausdorff dimension of exceptional times
and the exact noise sensitivity exponent for the hexagonal lattice
}\label{sect:GPS}

In \cite{GPS}, two of the main results (among many others)
were computing the exact noise sensitivity exponent for the crossing of
an $n\times n$ box
in the hexagonal lattice and the exact Hausdorff dimension of
the set of exceptional times for dynamical percolation on the 
hexagonal lattice. The latter number is the upper bound obtained in
Theorem \ref{th:SS}. 

\begin{thm} \label{th:SSbetter}\cite{GPS} 
For critical dynamical percolation on the hexagonal lattice, 
the Hausdorff dimension of the set of exceptional times is $31/36$.
In addition, for critical dynamical percolation on the square lattice,
there exist exceptional times of percolation.
\end{thm}

\begin{thm} \label{th:SSNoisebetter}\cite{GPS} 
Consider the sequence $\{f_n\}$ of indicator functions 
for a left to right crossing of an $n\times n$ square in the 
hexagonal lattice. Then 
$$
\STAB(\{f_n\})=\SENS(\{f_n\})=3/4.
$$
\end{thm}

\smallskip\noindent
(While we will not spell these out in detail,
for the square lattice, analogous results are obtained which relate
$\STAB(\{f_n\})$ and $\SENS(\{f_n\})$ with the expected number of pivotal
edges in large boxes.)

The value $31/36$ was the conjectured value in \cite{SS} and the
suggested value of $3/4$ was explained earlier.
The improvements here over Theorems \ref{th:SS} and \ref{th:SSNoise}
are due almost exclusively to providing much sharper results 
concerning the Fourier spectrum, both for crossing probabilities
and an annulus type event which is the relevant event for
studying exceptional times. These arguments do not for example use 
Theorem \ref{th:algorithm} or any variant of this.
This analysis is very long and intricate and so I will
only say a few words about it and even for that I will stick
to Theorem \ref{th:SSNoisebetter}.

\medskip\noindent
(Very vague) Outline of proof: \\
The much easier direction is to show that $\STAB(\{f_n\})  \le 3/4$. 
For a hexagon to be pivotal, there have to be,
starting from that hexagon, ``open paths''
to the left and right sides of the box and
``closed paths'' to the top and bottom sides of the box.
It is known that this has the same exponent as the 4-arm event. Looking at
points far from the boundary and using 
Theorem \ref{th:critexp}(ii), we obtain that $I(f_n)\ge n^{3/4+o(1)}$. 
In \cite{GPS}, it is shown that the boundary contributions can be 
controlled so that we indeed have $I(f_n)= n^{3/4+o(1)}$. 
Now Theorem \ref{th:influence} finishes the argument. 

The proof that $\SENS(\{f_n\}) \ge 3/4$ is significantly more difficult. 
By Theorem \ref{thm:meta}, we need to prove (\ref{eq:tightness}) and so
we need to obtain upper bounds on the lower tail of the distribution of
$|{\mathcal S}_n|$. Of course, we only care about the distribution of
$|{\mathcal S}_n|$ rather than the distribution of
${\mathcal S}_n$. However, it turns out that in order to study and analyze
the former, it is crucial to study the latter, which has much more 
structure and therefore more amenable to analysis. 
A first key step is for general Boolean functions and gives upper bounds on
\begin{equation}\label{eq:pivotalequation}
\P({\mathcal S}\cap B\neq \emptyset={\mathcal S}\cap W),
\end{equation}
where $B$ and $W$ are disjoint subsets of the domain variables $1,\ldots,n_m$,
in terms of a more general notion of pivotality. While one needs such a result
for all $W$, looking at the two special cases 
where $W$ is $\emptyset$ or $B^c$ illustrates well this relationship.
This general result gives in 
the context of percolation that if $B$ is a connected set of hexagons,
$\P({\mathcal S}\cap B\neq \emptyset)$ is at most 4 times the probability
of having 4 alternating arms from the set $B$ out to the boundary of our 
$n\times n$ box and 
$\P(\emptyset\neq {\mathcal S}\subseteq B)$ is at most 4 times the 
previous probability squared. This starts to get the ball rolling
as it relates the difficult spectral measure to things that are a
little bit more concrete.

Now let $g_r:=r^2\alpha_4(r)$ which is close to the expected number of
pivotals for a left to right crossing in an $r\times r$ box or equivalently
the expected size of the spectral sample for this event. This grows
to $\infty$ with $r$. One shows that
\begin{equation}\label{eq:tightnessspecialcase}
\P(|{\mathcal S}_n| < g_r,{\mathcal S}_n\neq \emptyset)
\asymp (n/r)^2 (\alpha_4(n)/\alpha_4(r))^2.
\end{equation} 
It is not hard to show, using the fact that the 
4 arm exponent is $5/4$, that (\ref{eq:tightness}) holds and that
one can take $\delta$ to be $\gamma^{3/2+\epsilon}$ for any fixed $\epsilon>0$.
While we want the upper bound, it is instructive
to see how the lower bound is obtained which is as follows. We break
the $n\times n$ square into about $(n/r)^2$ $r\times r$ squares.
It turns out that in the percolation context and with $B$ being
an $r\times r$ square and $W=B^c$, the upper bound 
on (\ref{eq:pivotalequation}) is shown to also be a lower bound 
(up to constants) and so for each $r\times r$ square $B$,  
the probability that the spectrum is nonempty
and sits inside $B$ can been shown to be at least 
$\Omega(1)(\alpha_4(n)/\alpha_4(r))^2$.
Next it is shown that, conditioned on the spectrum intersecting
such a box and no others, there is a uniform lower bound on the
probability that the spectral size within that box is at most $O(1)g_r$.  
Since, as we vary the $r\times r$ square, 
we obtain $(n/r)^2$ disjoint events, we obtain the (much easier)
lower bound of (\ref{eq:tightnessspecialcase}).

For the upper bound (which is much harder), we again break the 
$n\times n$ square as above and look at the number of $r\times r$ squares 
which intersect the spectral sample. Call this number $X_{n,r}$.
Using a very difficult geometric induction argument, it is shown that
$$
\P(X_{n,r}=k)\le g(k) (n/r)^2(\alpha_4(n)/\alpha_4(r))^2
$$
where $g(k)$ grows slower than exponentially (but faster than any 
polynomial). It turns out that there is a ``converse'' of what we wrote
in the previous paragraph which is that
conditioned on the spectrum touching a box, there is uniform lower bound 
on the probability that the size of the spectrum is at least
$\Omega(1)g_r$. This involves quite difficult percolation arguments.
Given $X_{n,r}=k$, if it were the case that the sizes of the
spectrum in the $k$ different $r\times r$ boxes which the spectrum hits
were independent, then the probability that all of them have size 
less than $\Omega(1)g_r$ would be
at most $c^k$ for some $c<1$. Since $\sum_k g(k)c^k <\infty$, we would be done.
The problem is, under this conditioning, the spectrum in the different boxes 
are not independent and have a very complicated dependency structure.
It is however shown that, while it is difficult to deal with the spectrum
in one box conditioned on its behavior elsewhere, it is possible to deal
with the spectrum in one box conditioned on it hitting that box and
not intersecting some other fixed set. This together with a novel 
type of large deviations argument allows us to carry out the upper 
bound in (\ref{eq:tightnessspecialcase}).
$\QED$

In proving the above, other interesting results
are obtained. For example, it is shown for percolation on $\Z^2$ 
that rerandomizing a small portion of only the vertical edges 
is sufficient for the correlation to go to 0. This result
suggests new interesting directions concerning
noise sensitivity for Boolean functions when only
some of the bits are rerandomized.

\section {Sensitivity of the infinite cluster
in critical percolation}\label{sect:PSS}

Except in the second part of Theorem \ref{th:G1G2},
all results in this paper so far dealt with graphs which
do not percolate at criticality.
It turns out that if we deal with graphs which do
percolate at criticality and ask if there are
exceptional times of {\em nonpercolation}, the structure
of the problem is quite different. In addition,
it seems to me that the set of exceptional times
in this case might have similar structure to the set
of ``fast points'' for Brownian motion; see the discussion
at the end of Section \ref{sect:tree}.

In \cite{PSS}, among other things,
a fairly detailed study of this
question was made for spherically symmetric trees. 
Special cases of the two main results of that paper
are the following. In this section (only), we are
dropping the assumption of homogeneous edge probabilities but will
assume all edge probabilities are bounded away from 0 and 1. The definition
of $w_n$ from Section \ref{sect:tree} should be modified in the obvious way.

\begin{thm} \label{th:sstrees1}
Consider a spherically symmetric tree with spherically symmetric edge 
probabilities (meaning all edges at a given level have the same
retention probability). Let $w_n$ be as in Theorem \ref{th:lyons}.

\smallskip\noindent
(i) If
$$
\lim_n \frac{\ew n}{n(\log n)^{\alpha}}=\infty
$$
for some $\alpha > 2$, then there are no exceptional times of  nonpercolation.

\smallskip\noindent
(ii) If
$$
{\ew n}\asymp {n(\log n)^{\alpha}}
$$
for some $1<\alpha \le 2$, then there are exceptional times of nonpercolation.

Note that in both of these regimes, Theorem \ref{th:lyons}
tells us that there is percolation at a fixed time.
\end{thm}

\medskip\noindent
{\em Remarks:} \\
(1) To see a concrete example, if
we have a tree with $|T_n| \asymp 2^nn(\log n)^\alpha$ and
$p=1/2$ for all edges, then if $\alpha>2$, we are in case (i)
while if $\alpha\le 2$, we are in case (ii). (Note Lyons'
theorem tells us that $p_c=1/2$ in these cases.) \\
(2) The theorem implies that if $\ew n\asymp n^\alpha$ with $\alpha > 1$, 
then there are no exceptional times of nonpercolation, while note that
if $\ew n \asymp n$, then Theorem \ref{th:lyons}
implies that there is no percolation
at a fixed time.
Hence, if we only look at the case where $w_n\asymp  n^\alpha$
for some $\alpha\ge 1$, we do not see the dichotomy (within
the regime where we percolate at criticality) that we are after
but rather we see the much more abrupt transition from not percolating at
a fixed time to percolating at all times.
Rather, Theorem \ref{th:sstrees1} tells us that
we need to look at a ``finer logarithmic scale'' to see this
``phase transition'' of changing from percolating at a fixed time but
having exceptional times (of nonpercolation) to percolating at all times.

\medskip
Interestingly, it turns out that even within the regime where
there are no exceptional times of nonpercolation,
there are still two very distinct dynamical behaviors 
of the process, yielding another phase transition.

\begin{thm}\label{th:dyntrans}
Consider a spherically symmetric tree $T$ of bounded degree and let $w_n$ 
be as in the previous result.

\smallskip\noindent
(i) When $\sum_{n=1}^\infty n\,\ew n^{-1}<\infty$, a.s.\ the set of
times $t\in[0,1]$ at which the root percolates
has finitely many connected components. (This holds for example
if $\ew n\asymp n^\theta$ with $\theta> 2$ as well as for supercritical
percolation on a homogeneous tree.)

\smallskip\noindent
(ii) If $\ew n\asymp n^\theta$, where
$1<\theta<2$, then with positive probability the set of times $t\in[0,1]$
at which the root percolates has infinitely many connected components.
The same occurs if $\ew n\asymp n(\log n)^\alpha$ for any $\alpha > 1$.
\end{thm}

\medskip\noindent
{\em Remarks:}\\
(1) If $\ew n\asymp n^2$, we do not know the 
answer. A first moment calculation suggests that there should be 
infinitely many connected components with positive probability but
the needed inequality for a successful second moment argument fails. \\
(2) It is easy to show that for any graph, if there 
are exceptional times of nonpercolation,
then the set of times $t\in[0,1]$ at which a 
fixed vertex percolates is totally disconnected and 
hence has infinitely many connected components
with positive probability. \\
(3) We will not indicate here any proofs but we will mention one word
about Theorem \ref{th:dyntrans} since the critical exponent of 2 
there derives
from a difference in the ordinary model in these two regimes. 
Namely, the expected number of pivotal edges for the 
event that the root percolates is infinite for $\theta \le 2$
but finite for $\theta > 2$.

\section {Dynamical percolation and the incipient infinite cluster}
\label{sect:IIC}

\subsection {The incipient infinite cluster}

We know that on $\Z^2$ and on the hexagonal lattice, there is no
infinite cluster at $p_c$. Nonetheless, physicists and others have tried to
talk about the ``infinite cluster on $\Z^2$ 
containing the origin at criticality''.
This was made sense of by Kesten in \cite{Keincipient} where the following
result was obtained. $\Lambda_n$ is the box of size $n$ centered at the origin.

\begin{thm}\label{th:incipient} \cite{Keincipient} The limiting measures 
$$
\lim_{p\downarrow 1/2} \P_p(\cdot\mid 0\leftrightarrow \infty)
$$
and
$$
\lim_{n\to\infty} \P_{1/2}(\cdot\mid 0\leftrightarrow \partial\Lambda_n)
$$
both exist and are equal. 
\end{thm}

This limiting measure is referred to as the {\bf incipient infinite cluster}.
Properties of this were also obtained and furthermore, in
\cite{Kerw}, Kesten showed that random walk on the incipient infinite cluster
is subdiffusive.

\subsection {The incipient infinite cluster and dynamical percolation}

It was asked quite early on whether the configuration for dynamical
percolation at a (properly chosen) 
exceptional time at which the origin percolates
(assuming there are such exceptional times) 
should have the same distribution as the incipient infinite
cluster of the previous subsection. 
See the discussion concerning this question in
\cite{Haggsurvey}. This question was answered in a very satisfactory manner
by Hammond, Pete and Schramm in \cite{HPSincip}, a paper which is
in the process of being written up. We sketch here a part of what was
done in this paper.

The first key step in being able to ``find'' the incipient infinite cluster 
inside of dynamical percolation is to be able to define a 
{\em local time} for when the origin is percolating. There are two different
approaches used to define a local time.

The first approach is as follows.
Let $A_{R,t}$ be the event that at time $t$ there is an open path 
from the origin to distance $R$ away and 
let $\calT_R$ be the random set of times at which 
$A_{R,t}$ occurs. Define a random measure $\mu_R$ on $\R$ by
$$
\mu_R:=\frac{1}{\alpha_1(R)}  \calL_{\calT_R}
$$
where $\calL_F$ refers to Lebesgue measure restricted to the set $F$
and $\alpha_1(R)=\P(A_{R,t})$. Clearly, $\E[\mu_R([a,b])]=b-a$.
We then let $R$ tend to infinity.

The second approach (which turns out to be equivalent)
is as follows and is closer in spirit to 
Kesten's original definition of the incipient infinite cluster.
Consider ordinary percolation and 
let $S$ be a collection of hexagons and let $\omega^S$ be the
percolation realization restricted to $S$. We want to measure
in some sense how much $\omega^S$ ``helps percolation to 
infinity''. This is made precise by the limit 
$$
\lim_{R\to\infty}\frac{\P(A_R\mid  \omega^S)}{\P(A_R)}
$$
which is easily shown to exist using Theorem \ref{th:incipient}.

Calling this limit $f(\omega^S)$, let $M_r(\omega):=f(\omega^{B_r})$
where $B_r$ is the ball of radius $r$ around the origin.
Finally let 
$$
\nu_r([a,b]):=\int_a^b M_r(\omega_s) ds.
$$

\begin{thm}\label{th:local} \cite{HPSincip} \\
(i) For all $ a< b$, $\mu_R([a,b])$ converges a.s.\ and in $L^2$ to
a limit as $R$ goes to $\infty$, which we call $\mu_\infty([a,b])$. \\
(ii) For all $ a< b$, $\nu_r([a,b])$ converges a.s.\ and in $L^2$ to
a limit as $r$ goes to $\infty$, which we call $\nu_\infty([a,b])$. \\
(iii) The two above limits are a.s.\ the same.
\end{thm}

\medskip\noindent
Clearly the limiting measure $\mu_\infty$ is supported on
the set of exceptional times at which the origin percolates
(the latter known to be a nonempty closed set).
It is not known if the support of the measure is exactly the set
of exceptional times. It is explained that $\mu_R([a,b])$ is a martingale
which implies its a.s.\ convergence. Using estimates from \cite{SS},
one can see it is also $L^2$ bounded which gives the $L^2$ convergence.
The convergence in $L^2$ guarantees that the limiting measure is nonzero.
It is also shown that $\nu_r([a,b])$ is a martingale.

The final result links up the incipient infinite cluster
with dynamical percolation.

\begin{thm}\label{th:final} \cite{HPSincip}
Consider the random measure $\mu_\infty$ on $\R$ above and let $X$ 
be a Poisson process on $\R$ with ``intensity measure'' $\mu_\infty$.
Then the distribution of $\omega_0$ given $0\in X$ has the same 
distribution as the incipient infinite cluster.
\end{thm}

It is not so hard to make sense of the conditioning $0\in X$ even if
this event has probability 0; see Chapter 11 of \cite{Kal}. 
There are a number of other results in this paper which we do not detail 
here.

\section {The scaling limit of planar dynamical percolation}\label{sect:SL}

Before discussing the scaling limit of dynamical percolation,
it is necessary to first discuss the scaling limit of ordinary percolation.
There is a lot to be said here and this section will be somewhat
less precise than the earlier sections. Since even the formulations
can be quite technical, I will be, unlike in the previous sections, 
``cheating'' in various places.

Even before we state the scaling limit of percolation, we need to
first briefly explain the concept of conformal invariance and Cardy's
formula. Let $\Omega$ be a simply connected open domain in the plane and let
$A,B,C$ and $D$ be 4 points on the boundary of $\Omega$ in clockwise order.
Scale a 2-dimensional lattice, such as $\Z^2$ or the
hexagonal lattice, by $1/n$ and perform critical percolation on this 
scaled lattice.
Let $\P(\Omega, A,B,C,D,n)$ denote the probability that, in the
$1/n$ scaled hexagonal lattice, there is a white path of hexagons inside
$\Omega$ going from the boundary of $\Omega$ between $A$ and $B$
to the boundary of $\Omega$ between $C$ and $D$. The first half
of the following conjecture
was stated in \cite{LPS} and attributed to Michael Aizenman. The second half
of the conjecture is due to Cardy \cite{Cardy}.

\begin{conj}
(i) For all $\Omega$ and $A,B,C$ and $D$ as above,
$$
\P(\Omega, A,B,C,D,\infty):=\lim_{n\to\infty}\P(\Omega, A,B,C,D,n)
$$
exists and is conformally invariant in the sense that
if $f$ is a conformal mapping, then
$\P(\Omega, A,B,C,D,\infty)= \P(f(\Omega), f(A),f(B),f(C),f(D),\infty)$. \\
(ii) There is an explicit formula (not stated here)
for $\P(\Omega, A,B,C,D,\infty)$, called Cardy's formula,
when $\Omega$ is a rectangle and $A,B,C$ and $D$ are the 4 corner points.
(Since every $\Omega$, $A,B,C$ and $D$ can be mapped to a unique such 
rectangle (with $A,B,C$ and $D$ going to the 4 corner points), 
this would specify the above limit in general assuming conformal invariance.)
\end{conj}

Cardy's formula is quite complicated involving hypergeometric functions but
Lennart Carleson realized that assuming conformal invariance, there is
a nicer set of ``representing'' domains with four specified points for
which the limit has a much simpler form.
Namely, if $\Omega$ is an equilateral 
triangle (with side lengths 1), $A,B$ and $C$ the three corner points
and $D$ (on the line between $C$ and $A$)
having distance $x$ from $C$, then the above probability would just be
$x$.  Using Carleson's reformulation of Cardy's formula, Smirnov proved 
the above conjecture for the hexagonal lattice.

\begin{thm} \label{th:smirnov}\cite{Smirnov}
For the hexagonal lattice, both (i) and (ii) of the above conjecture are
true.
\end{thm}

This conjecture is also believed to hold on $\Z^2$ but is not (yet) 
proved in that case. In \cite{Schramm},
Schramm described what the interfaces between whites and blacks should be
as the lattice spacing goes to 0, assuming conformal invariance. In the
appropriate formulation, it should be an SLE$_6$ curve. 
Smirnov \cite{Smirnov}
proved this convergence for one interface and Camia and Newman
\cite{CN} proved a ``full scaling limit'', which is a description of the
behavior of all the interfaces together. The critical exponents described
in Section \ref{sect:crit} are proved by exploiting the SLE description
of the interfaces.  All of the above is described well in \cite{W}.

It turns out, in order to obtain a scaling limit of dynamical percolation,
a different description, due to Schramm and Smirnov,
of the scaling limit of ordinary percolation is preferable. 
A quad $Q$
is a subset of the plane homeomorphic to a disk together with its boundary
partitioned into 4 continuous pieces. A configuration for the scaling limit
is described by a family of $0$-$1$ random variables $X_Q$ indexed by
the quads $Q$ where $X_Q=1$ means there is a crossing from the first to
the third boundary piece of $Q$ using white hexagons. 
Equivalently, an element of the state space is 
a collection of quads (satisfying a certain necessary monotonicity 
condition and suitably measurable) where the collection of quads 
represents the quads which are ``crossed''. An important advantage of
this state space, which we denote by ${\mathcal S}$, is that it is
a compact metrizable space. This avoids the need for tightness arguments.
If we perform percolation on the
$1/n$ scaled hexagonal lattice, the $X_Q$'s are well-defined random
variables and so we obtain a probability measure $\mu_n$ on
${\mathcal S}$. It turns out that the sequence of probability measures 
$\mu_n$ have a limit $\mu_\infty$, which we call the scaling limit.
Note that for each $n$, there is essentially a 1-1 correspondence between
percolation realizations on the $1/n$ scaled lattice and
elements in ${\mathcal S}$. In the limit however, 
we can not talk anymore about percolation configurations but rather the 
only information left is which quads are ``crossed''. (The latter is sort 
of the ``macroscopic'' information.)

We now move to a scaling limit for dynamical percolation. Consider 
dynamical percolation on the $1/n$ scaled hexagonal lattice.
We can think of our dynamical percolation
on this scaled lattice as a process $\{\eta^n(t)\}_{t\in\R}$ taking
values in ${\mathcal S}$. We can try to let $n$ go to infinity, in which
case, the marginal for each $t$ will certainly go to $\mu_\infty$.
However due to the noise sensitivity of percolation, for any fixed $s <t$
and any quad $Q$, the crossing of $Q$ at time $s$ and at time $t$ will become
asymptotically independent as $n\to\infty$ which implies that the processes
$\{\eta^n(t)\}_{t\in\R}$ converge to something which is 
``independent copies of $\mu_\infty$ for all different times''. This
is clearly not what we want. In order to obtain a nontrivial limit,
we should slow down time by a factor of $1/(n^2\alpha_4(n))$ 
where $\alpha_4(n)$ is explained in Section \ref{sect:GPS}.

\begin{thm} \label{th:scalinglimit}\cite{GPS2}
Let $\sigma^n_t:=\eta^n_{tn^{-2}\alpha^{-1}_4(n)}$. Then 
$\{\sigma^n(t)\}_{t\in\R}$ converges in law to a process
$\{\sigma^\infty(t)\}_{t\in\R}$ under the topology of local uniform convergence.
Moreover, the limiting process is Markovian (which is not apriori obvious
at all).
\end{thm}

We explain now why the time scaling is as above. Consider the quad 
corresponding to the unit square together with its four sides. The expected 
number of pivotals for a left to right crossing of this quad on the 
$1/n$ scaled lattice is $\asymp n^2\alpha_4(n)$. 
Next, the above time scaling updates each hexagon in unit time with 
probability about $1/n^2\alpha_4(n)$ and therefore, by the above, 
updates on average order 1 pivotals. It follows from the main result in
\cite{GPS} that if we scale
``faster'' than this, you get a limit which is independent
at different times and if we scale ``slower'' than this, you get a 
limit which is constant in time. This scaling is such that the 
correlation between this event occuring at time 0 and occuring at time 
1 stays strictly between 0 and 1, which is something we of course want.

In the above paper, there is another model which is studied called 
near-critical percolation for which a scaling limit is proved in a similar
fashion. Near-critical percolation can be thought of as
a version of dynamical percolation, where sites are 
flipping only in one direction and hence after time $t$, 
the density of open sites is about $1/2 + t/(n^2\alpha_4(n))$.

To prove Theorem \ref{th:scalinglimit}, there are two key steps. 
The first key step is a stability type result. In very vague terms, it says 
that if we look at the configuration at time 0, ignore the evolution of
the hexagons which are initially pivotal only for crossings which occur at
macroscopic scale at most $\rho$ (which should be thought of as 
much larger than the lattice spacing) but observe the evolution of
the hexagons which are initially pivotal for crossings which occur at
macroscopic scale larger than $\rho$,
then we can still predict quite well 
(i.e., with high probability if $\rho$ is small)
how the crossings evolve on macroscopic scale 1.

Since we cannot see the individual hexagons which are initially pivotal
for crossings at scale at least $\rho$
in the scaling limit, in order for this stability result to be useful, we
need at least that the number of such hexagons in a given region can be
read off from the scaling limit, so that we know the rate with which
macroscopic crossings are changing. This is obtained by
the second key step. For every $n$, look at ordinary percolation
on the $1/n$ scaled lattice and consider the set of hexagons
from which emanate 4 paths of alterating colors to macroscopic-distance 
$\rho$ away. (These are the hexagons which, if flipped, change the state
of a crossing at scale $\rho$.) Let 
$\nu_n^\rho$ be counting measure on such points divided by $n^2\alpha_4(n)$.
This scaling is such that the expected measure of the unit square is
of order 1 (as $n \to\infty$ if $\rho$ is fixed).

\begin{thm} \label{th:measure}\cite{GPS2}
For all $\rho$, $(\omega_n,\nu_n^\rho)$ converges,
as $n \to\infty$, to a limit $(\omega_\infty,\nu_\infty^\rho)$ {\it and}
$\nu_\infty^\rho$ is a measurable function of $\omega_\infty$. 
\end{thm}

The argument of this result is very difficult. I say however a word on how 
this gives a nice construction of the scaling limit of dynamical percolation.
The idea is that one constructs a dynamical percolation realization
by first taking a realization from the $t=0$ scaling
limit (or equivalently $\omega_\infty$ above), looks at 
$\nu_\infty^\rho$ (which requires no further randomness because 
$\nu_\infty^\rho$ has been proved to be a function of $\omega_\infty$)
and then builds an appropriate Poisson point process over the random measure 
$\nu_\infty^\rho\times dt$ which will be used to 
dictate when the different crossings of 
quads of scale at least $\rho$ will change their state. This is done for
all $\rho$ and the stability result is needed to ensure that the 
process described in this way is in fact the scaling limit of dynamical 
percolation. The idea that dynamical percolation 
(and near-critical percolation) could possibly be built up in this manner
from a Poisson process over a random measure was suggested by
Camia, Fontes and Newman in the paper \cite{CFN}.

It is also shown in \cite{GPS2} that a  type of ``conformal covariance''
for the above measures $\nu_\infty^\rho$ holds but this will not be detailed 
here. The argument for Theorem \ref{th:measure} also proves the existence of 
natural time-parametrizations for the SLE$_6$ and SLE$_{8/3}$ curves, a 
question studied in \cite{L} for general SLE$_\kappa$.

\section {Dynamical percolation for interacting particle systems}
\label{sect:BS}

In \cite{BS}, certain interacting particle systems 
were studied and the question of whether
there are exceptional times at which the percolation
structure looks different from that at a fixed time was
investigated. The two systems studied were the
contact process and the Glauber dynamics for the Ising model.
Most of the results dealt
with noncritical cases but since the dynamics
are not independent, a good deal more work was needed
compared to the easy proof of Proposition \ref{pr:noncrit}. 

However, there was one very interesting case
left open in this work which was the following.
Consider the Glauber dynamics for the critical
$2$-dimensional Ising model. First note that
$2$ dimensions is special for the Ising model
compared to higher dimensions since it is known that for $\Z^2$
the critical value for phase transition is the same as the critical
value for when percolation occurs. The Ising model
does not percolate in $2$ dimensions at the
critical value and in view of Theorem \ref{th:SS},
it is natural to ask if there are exceptional
times of percolation for the Glauber dynamics.
During a visit to Microsoft, G\'{a}bor Pete 
suggested to me that this might not be the case
since the scaling limit of the Ising model was
conjectured to be SLE$_3$, which unlike 
SLE$_6$ (the scaling limit of percolation), does not
hit itself. This yields that there should not be
so many pivotal sites which suggests no exceptional times.
A few months later in Park City, Christophe Garban 
sketched out a ``1st moment argument'' for me which would
give no exceptional times under certain assumptions.

In any case, whether or not one can actually prove this,
it now seems clear that there are no exceptional times
for the Glauber dynamics. This seems to be
a beautiful example of how the
qualitative difference in the behavior of the two
respective SLE scaling limits (self-intersections versus 
not) yields a fundamental difference between the natural 
dynamics on ordinary percolation and the natural dynamics on the Ising model
concerning whether there are exceptional times of percolation or not.

\subsection*{Acknowledgments}
I would like to thank Erik Broman, Christophe Garban,
Olle H\"aggstr\"om, Alan Hammond, Russell Lyons and G\'{a}bor Pete for
very many helpful comments on earlier versions of the manuscript.
I appreciate in particular Christophe Garban, Alan Hammond and G\'{a}bor 
Pete explaining to me their unpublished work. A part of this work was
done while visiting the Mittag-Leffler Institute in Djursholm, Sweden,
the same institute where dynamical percolation was started.

\end{document}